\documentclass[12pt]{amsart}

\theoremstyle{remark}
\newtheorem*{exercise}{Exercise for the reader}

\usepackage{url}

\usepackage{tikz}
\usetikzlibrary{decorations.pathmorphing, decorations.markings}

\begin{document}

\title{Mocposite functions}
\author{Harold P. Boas}
\address{Department of Mathematics \\ Texas A\&M University \\ College Station, TX 77843 \\ USA}
\email{boas@tamu.edu}

\maketitle

An engineering student and a mathematics student walk into a bar. Instead of carding the students, the bartender offers them  free drinks for a correct answer to the question, ``Is the function \(\sqrt{1-z^{2}}\) even or odd?'' The engineering student shouts out, ``Even, of course!'' Noticing the bartender's sphinx tattoo, the mathematics student slyly says, ``My answer is: Yes.'' 
Although the smart-aleck second answer arguably is less wrong than the first, 
the bartender throws both students into the street and orders them to stay away until they have studied analytic continuation. 

The surprise is that \(\sqrt{1-z^{2}}\) appears in some applications as an odd function! This statement seems absurd at first sight, for \(1-z^{2}\) is manifestly even, and composing an even function with any subsequent operation preserves evenness. 
The startling resolution of this paradox is that \(\sqrt{1-z^{2}}\) only pretends to be a composite function but actually is not one.
I propose that such mock composite functions be called mocposite. 

This article studies \(\sqrt{1-z^{2}}\) as a means of entering the looking-glass world \cite{mirror} of mocposite functions, where even is odd, and odd is ordinary. Some prior acquaintance with the elements of complex analysis will make the reader's passage smoother. 
My tale includes both a caution on confusing conventions and a pedagogical praise of pedantry. 

\section{Indices and surds}
Understanding \(\sqrt{1-z^{2}}\) requires first coming to grips with the notation for square roots. 
The peculiar symbol~\(\surd\) dates from sixteenth-century Germany, according to Florian Cajori \cite[\S\S316--338]{cajori}, and the juxtaposition of the horizontal grouping bar (vinculum) is a subsequent innovation of Ren\'e Descartes---one of his most enduring and most regrettable contributions to mathematics. Why not use exponent~\(1/2\) to denote a square root? The exponential form is both cleaner than~\(\surd\) to typeset and consistent with the standard notation  for other powers.

A thornier problem than the notation is the ambiguity inherent in the concept of square root, for every number has two square roots. If you object that the number~\(0\) is an exception having a single square root, then observe that what \(\sqrt{z}\) really means is a solution~\(w\) to a particular quadratic equation: namely, \(w^2-z=0\). Every quadratic equation has two solutions, counting multiplicity. 

Nonetheless, there is one case in which everybody agrees that the symbol~\(\sqrt{z}\) has a unique meaning. When \(z\)~happens to be a positive real number, convention dictates that \(\sqrt{z}\) always denotes the positive square root of~\(z\). But if \(z\)~is a negative real number (or, worse, a nonreal number), then confusion can and does arise. 

A quantity~\(i\) whose square equals~\(-1\) is fundamental to complex analysis, so neither the existence nor the uniqueness of~\(i\) should pass without comment. In the influential terminology of Descartes~\cite[p.~380]{descartes}, nonreal solutions of polynomial equations are ``imaginary'' in the sense of existing only in the imagination. The device of giving imaginary numbers a concrete existence as ordered pairs of real numbers (equipped with a suitable multiplication) is due to William Rowan Hamilton~\cite{hamilton} two hundred years after Descartes. The imaginary unit~\(i\) has an alternative realization, invented by Augustin-Louis Cauchy~\cite{cauchy}, that can be expressed in modern language as the equivalence class of the indeterminate~\(x\) in the algebraic structure \(\mathbb{R}[x]/(x^{2}+1)\), the quotient of the ring of polynomials with real coefficients by the ideal consisting of polynomials that have \(x^{2}+1\) as a factor. 

Authors who wish to have the letter~\(i\) available as a summation index often write a complex variable in the form \(x+y\sqrt{-1}\) instead of \(x+yi\), innocently imagining (I suppose) that the symbol~\(\sqrt{-1}\) has a unique meaning rather than two possible values. An inevitable consequence of this belief would be that \(\sqrt{-4}\) has a unique meaning (namely, \(2\sqrt{-1}\,\)), and more generally that \(\sqrt{z}\) is well defined for~\(z\) everywhere on the negative part of the real axis. Since this set is precisely the standard branch cut across which the complex square-root function is discontinuous, such authors are implicitly constructing an edifice on top of a fault line and hoping that no earthquake occurs. 

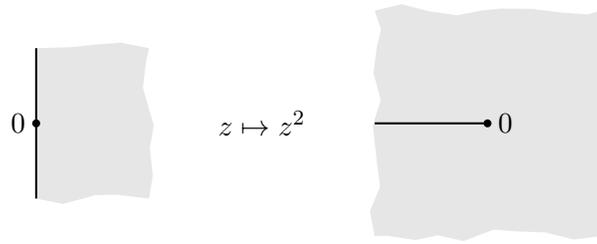
\begin{figure}
\begin{center}
\small
\begin{tikzpicture}[thick]
\fill[gray!20, decoration=random steps] decorate {(0,-1) -- (1.5,-1) -- (1.5,1) -- (0,1)} -- cycle;
\fill (0,0) circle (1.5pt);
\draw (3,0) node {$z\mapsto z^{2}$};
\draw (0,-1) -- (0,1);
\draw (0,0) node[anchor=east] {$0$};
\begin{scope}[xshift=6cm]
\fill[gray!20, decoration=random steps, decorate] (-1.5,-1.5) rectangle (1.5,1.5); \draw (0,0) -- (-1.5,0); 
\fill (0,0) circle (1.5pt);
\draw (0,0) node[anchor=west] {$0$};
\end{scope}
\end{tikzpicture}
\end{center}
\caption{The squaring function}
\label{fig:squaring}
\end{figure}

The standard square-root function arises by considering an inverse of the function that sends a complex number~\(z\) to the image~\(z^{2}\). As indicated in Fig.~\ref{fig:squaring}, this squaring function maps the open right-hand half-plane (where the real part of~\(z\) is positive) bijectively onto the complex plane with a left-hand slit along the real axis from~\(0\) to~\(-\infty\). The ``principal branch of~\(\sqrt{z}\,\)'' means the inverse of this squaring function. Being the inverse of a holomorphic (that is, complex-analytic) function, the principal branch of~\(\sqrt{z}\) is a holomorphic function too. More generally, ``a branch of~\(\sqrt{z}\,\)'' means a holomorphic function~\(f\) such that \((f(z))^{2}=z\) for every~\(z\) in some prescribed domain in the complex plane. 

Not every domain supports a branch of~\(\sqrt{z}\). The obstruction is the existence in the domain of a simple closed curve~\(\gamma\) that surrounds the origin. Indeed, if \((f(z))^{2}=z\) for every~\(z\) in some domain, then the chain rule implies that \(2ff'=1\), so 
\begin{equation*}
\frac{f'(z)}{f(z)} = \frac{1}{2(f(z))^{2}} = \frac{1}{2z}.
\end{equation*}
If \(C\) denotes the image of~\(\gamma\) under~\(f\), then 
\begin{equation*}
\frac{1}{2\pi i} \int_{\gamma} \frac{f'(z)}{f(z)}\,dz = \frac{1}{2\pi i} \int_{C} \frac{1}{w}\,dw,
\end{equation*}
which equals the winding number of the curve~\(C\) about~\(0\): namely, a particular integer. On the other hand, this integer equals
\begin{equation*}
\frac{1}{2\pi i} \int_{\gamma} \frac{1}{2z}\,dz,
\end{equation*}
which is half the winding number of~\(\gamma\) about~\(0\). Hence the existence of~\(f\) precludes the existence of a curve~\(\gamma\) in the domain with winding number~\(1\) about~\(0\), since \(1/2\) is not an integer. 

When \(g\) is a holomorphic function, ``a branch of \(\sqrt{g}\,\)'' means a holomorphic function~\(f\) such that \((f(z))^2 = g(z)\) for every~\(z\) in some prescribed domain. A subtle but crucial point is that the existence of a branch of~\(\sqrt{g}\) does not necessarily entail the existence of a branch of~\(\sqrt{z}\) on the image of~\(g\). If \(g(z)=z^2\), for instance, then there is a branch of~\(\sqrt{g}\) on the entire complex plane (namely, the identity function), but there is no branch of~\(\sqrt{z}\) on the image of~\(g\) (for that image is the entire complex plane). 

If a domain supports a branch of \(\sqrt{z}\), then the negative of that function is another branch. Consequently, the value of the expression~\(\sqrt{z}\) when \(z=4\) is not necessarily equal to~\(\sqrt{4}\) (since \(\sqrt{4}\) conventionally is positive). Do you sense the other-worldly weirdness wafting from the standard notation for square roots?

Pedants distinguish between the name of a function, say \(\cos\), and the value of a function at a point, say \(\cos z\). Most authors, however, use the notation \(\cos z\) ambiguously to mean either ``the value of the function \(\cos\) at the point~\(z\)'' or ``the function that sends the variable~\(z\) to the value \(\cos z\).'' My father was fond of pointing out that the second usage corresponds to a standard trope of classical rhetoric: synecdoche is the figure of speech in which a part stands for the whole. Normally, no confusion arises from naming a function by a generic value, but \(\sqrt{1-z^{2}}\) presents a dramatic exception, as I shall demonstrate now. 

\section{And this was odd}
\label{odd}
When I was an undergraduate, back in the days when the distinguished mathematician John Tate had won only the first of his many major awards, I heard him declare that ``2~is an odd prime'' (an entirely reasonable statement in the context of algebraic number theory). I intend to make the case that \(\sqrt{1-z^{2}}\) is an odd function (in both senses of the word ``odd''). 

Keep in mind that what the expression \(\sqrt{1-z^{2}}\) means is a function~\(f\) such that \( (f(z))^{2} = 1-z^{2}\) for all~\(z\) in some specified domain. Introducing the variable~\(w\) to represent \(f(z)\) converts the equation into the following form:
\(w^{2}+z^{2}=1\).
This relation defines a certain subset of the space~\(\mathbb{C}^{2}\) of two complex variables, a subset that some readers may wish to think of as a Riemann surface (a one-dimensional complex manifold). The implicit-function theorem implies that \(w\)~can be expressed as a holomorphic function of~\(z\) locally near each point on the surface at which \(w\)~is different from~\(0\) (equivalently, \(z\)~is different from~\(\pm 1\)). 

Since the equation is symmetric with respect to the two variables, there is no reason for~\(w\) to play a distinguished role. If the equation determines \(w\) as a function of~\(z\) in some region of the complex plane, then symmetry dictates that \(z\)~is the same function of~\(w\) in the identical region of the \(w\)-plane. Actually there must be two functions, negatives of each other, since the equation does not distinguish between \(w\) and~\(-w\) (or between \(z\) and~\(-z\)). 

Symmetry considerations thus give rise to the problem of prescribing a suitable subdomain of \(\mathbb{C} \setminus\{0, 1, -1\}\) and a bijective holomorphic function~\(f\) from that subdomain to itself such that \((f(z))^{2} = 1-z^{2}\) for every~\(z\). Moreover, the inverse function must be either the same function~\(f\) or its negative.

\begin{figure}
\begin{center}
\small
\begin{tikzpicture}[thick]
\fill[gray!20, decoration=random steps] decorate {(1.5,0) -- (1.5,1.5) -- (-1.5,1.5) -- (-1.5,0)} -- cycle;
\fill (0,0) circle (1.5pt);
\fill (1,0) circle (1.5pt);
\fill (-1,0) circle (1.5pt);
\draw (3,0) node {$z\mapsto 1-z^{2}$};
\draw (-1.5,0) -- (1.5,0);
\draw (0,0) node[anchor=north] {$0$};
\draw (1,0) node[anchor=north] {$1$};
\draw (-1,0) node[anchor=north] {$-1\phantom{-}$};
\begin{scope}[xshift=6cm]
\fill[gray!20, decoration=random steps, decorate] (-1.5,-2) rectangle (2,2);
\draw (-1.5,0) -- (1,0); 
\fill (1,0) circle (1.5pt);
\draw (1,0) node[anchor=north] {$1$};
\fill (0,0) circle (1.5pt);
\draw (0,0) node[anchor=north] {$0$};
\end{scope}
\end{tikzpicture}
\end{center}
\caption{Preparation for taking a square root}
\label{fig:pre}
\end{figure}
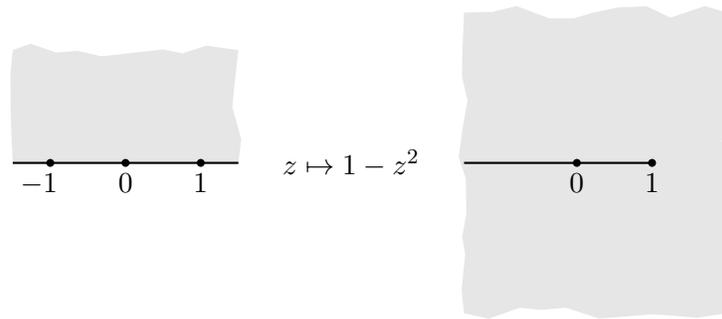

A reasonable initial step in the construction of~\(f\) is to define a branch of \(\sqrt{1-z^{2}}\) on the upper half-plane, the set where \(z\)~has positive imaginary part. Every complex number that is neither a positive real number nor~\(0\) is the square of exactly one such~\(z\). Therefore the function that sends~\(z\) to \(1-z^{2}\) maps the open upper half-plane bijectively onto the plane with a left-hand slit along the real axis from~\(1\) to~\(-\infty\). (See Fig.~\ref{fig:pre}.) This open region is a subset of the domain of the principal branch of the square-root function, so \(\sqrt{1-z^{2}}\) is well-defined on the upper half-plane as a composite function, say~\(f_{1}\). 

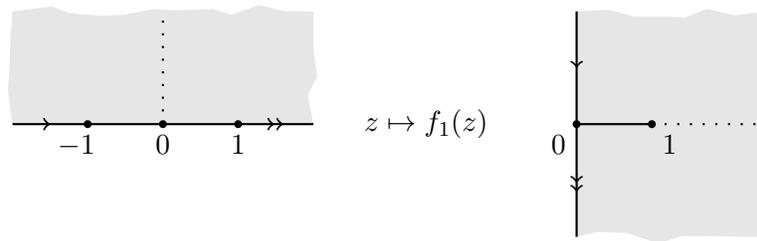
\begin{figure}
\begin{center}
\small
\begin{tikzpicture}[thick]
\fill[gray!20, decoration=random steps] decorate {(2,0) -- (2,1.5) -- (-2,1.5) -- (-2,0)} -- cycle;
\fill (0,0) circle (1.5pt);
\fill (1,0) circle (1.5pt);
\fill (-1,0) circle (1.5pt);
\draw (3.5,0) node {$z\mapsto f_{1}(z)$};
\draw[decoration={markings, mark=at position 0.5 with {\arrow{>}}}, postaction={decorate}] (-2,0) -- (-1,0);
\draw (-1,0) -- (1,0);
\draw[decoration={markings, mark=at position 0.6 with {\arrow{>>}}}, postaction={decorate}] (1,0) -- (2,0);
\draw (0,0) node[anchor=north] {$0$};
\draw (1,0) node[anchor=north] {$1$};
\draw (-1,0) node[anchor=north] {$-1\phantom{-}$};
\draw[loosely dotted] (0,0) -- (0,1.5);
\begin{scope}[xshift=5.5cm]
\fill[gray!20, decoration=random steps] decorate {(0,-1.5) -- (2.5,-1.5) -- (2.5, 1.5) -- (0,1.5)} -- cycle;
\draw[decoration={markings, mark=at position 0.5 with {\arrow{>}}}, postaction={decorate}] (0,1.5) -- (0,0);
\draw[decoration={markings, mark=at position 0.6 with {\arrow{>>}}}, postaction={decorate}] (0,0) -- (0,-1.5);
\draw (0,0) -- (1,0); 
\fill (1,0) circle (1.5pt);
\fill (0,0) circle (1.5pt);
\draw (0,0) node[anchor=north east] {$0$};
\draw (1,0) node[anchor=north west] {$1$};
\draw[loosely dotted] (1,0) -- (2.5,0);
\end{scope}
\end{tikzpicture}
\end{center}
\caption{$\sqrt{1-z^{2}}$ on the upper half-plane}
\label{fig:uhp}
\end{figure}

The image of~\(f_{1}\) is nearly identical to the image of the principal branch of the square root, except for removal of the image under the square-root function of the segment of the real axis from \(0\) to~\(1\). Since the square-root function maps that segment back to itself, the function~\(f_{1}\) maps the upper half-plane bijectively onto the right-hand half-plane with a slit along the real axis from \(0\) to~\(1\), as shown in Fig.~\ref{fig:uhp}. Notice that if \(y\)~is a positive real number, then \(f_{1}(iy) = \sqrt{1+y^{2}}\) (positive square root), so \(f_{1}\)~maps the part of the imaginary axis in the upper half-plane onto the unbounded interval of the real axis from \(1\) to~\(+\infty\). 

The next step---a nonunique process---is to extend the function~\(f_{1}\) to a larger domain. Here is one way to proceed. Observe that when \(z\)~is a point in the upper half-plane with real part greater than~\(1\) and imaginary part close to~\(0\), the point~\(z^{2}\) has the same properties. The point \(1-z^{2}\) then lies in the third quadrant close to the real axis. Taking the principal square root shows that the value \(f_{1}(z)\) lies in the fourth quadrant close to the imaginary axis. The upshot is that \(f_{1}\)~extends continuously to the unbounded interval of the real axis to the right of~\(1\) and maps this interval to the bottom half of the imaginary axis. Explicitly, the extension of~\(f_{1}\) maps an arbitrary real number~\(x\) greater than~\(1\) to the image \( -i\sqrt{x^{2}-1}\) (positive square root). Parallel reasoning shows that \(f_{1}\)~extends continuously to the unbounded interval of the real axis to the left of~\(-1\), and \(f_{1}\) maps this interval to the top half of the imaginary axis (Fig.~\ref{fig:uhp}). 

This situation admits application of the Schwarz reflection principle, the simplest method of analytic continuation discussed in a first course on complex analysis. The principle says that if a holomorphic function in the top half of a region symmetric with respect to the real axis extends continuously to an open subset of the real axis and takes real values there, then the function extends across that subset of the real axis to a function that is holomorphic in the whole symmetric region. Moreover, the extended function maps points that are symmetric with respect to the real axis to image points that are again symmetric with respect to the real axis.

Accordingly, the function~\(if_{1}\) extends by reflection to be holomorphic on the plane with a slit along the real axis from \(-1\) to~\(1\). Let \(f_{2}\) denote the corresponding extension of~\(f_{1}\) to this slit plane. Since the extension of~\(if_{1}\) maps pairs of complex-conjugate points to complex-conjugate image points, the function~\(f_{2}\) has the property that 
\begin{equation}
f_{2}(\overline{z}) = -\overline{f_{2}(z)}
\label{eq:sym}
\end{equation}
for every point~\(z\) in the slit plane. 

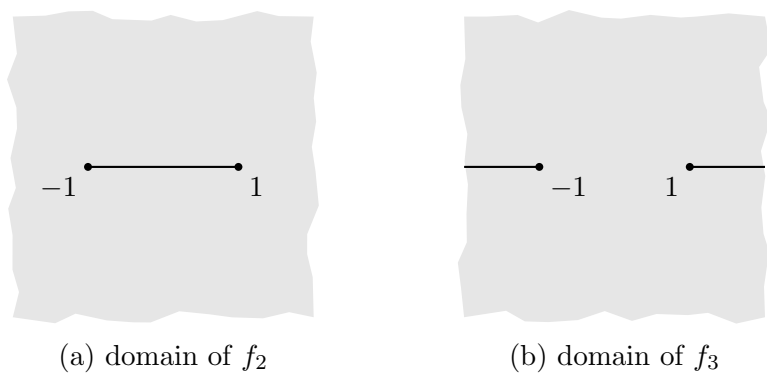
\begin{figure}
\begin{center}
\small
\begin{tikzpicture}[thick]
\fill[gray!20, decoration=random steps, decorate] (-2,-2) rectangle (2,2);
\draw (-1,0) -- (1,0);
\fill (-1,0) circle (1.5pt);
\fill (1,0) circle (1.5pt);
\draw (-1,0) node[anchor=north east] {$-1$};
\draw (1,0) node[anchor=north west] {$1$};
\draw (0,-2.2) node[anchor=north] {(a) domain of $f_{2}$};
\begin{scope}[xshift=6cm]
\fill[gray!20, decoration=random steps, decorate] (-2,-2) rectangle (2,2);
\draw (-2,0) -- (-1,0);
\draw (1,0) -- (2,0);
\fill (-1,0) circle (1.5pt);
\fill (1,0) circle (1.5pt);
\draw (-1,0) node[anchor=north west] {$-1$};
\draw (1,0) node[anchor=north east] {$1$};
\draw (0,-2.2) node[anchor=north] {(b) domain of $f_{3}$};
\end{scope}
\end{tikzpicture}
\end{center}
\caption{Two domains for $\sqrt{1-z^{2}}$}
\label{fig:domains}
\end{figure}

When \(z\) lies in the upper half-plane, \((f_{2}(z))^{2} = (f_{1}(z))^{2}=1-z^{2}\). Two holomorphic functions that agree on an open set agree identically on their common connected domain (by the identity principle from a first course on complex analysis), so \((f_{2}(z))^{2}=1-z^{2}\) on the whole plane with a slit along the real axis from \(-1\) to~\(1\). In other words, \(f_{2}(z)\) gives a well-defined meaning to \(\sqrt{1-z^{2}}\) on this slit plane, shown in Fig.~\ref{fig:domains}(a). 

What symmetry property does \(f_{2}\)~have?
Letting \(y\)~be a positive real number and setting \(z\) equal to~\(iy\) in equation~\eqref{eq:sym} shows that \(f_{2}(-iy) = -\overline{f_{2}(iy)}\). Since \(f_{1}\) (hence~\(f_{2}\)) takes real values on the top half of the imaginary axis, the preceding equation implies that \(f_{2}(-iy)=-f_{2}(iy)\). In other words, the expression \(f_{2}(-z)+f_{2}(z)\) is identically equal to zero when \(z\)~lies on the top half of the imaginary axis. Since zeroes of nonconstant holomorphic functions are isolated, the sum \(f_{2}(-z)+f_{2}(z)\) is identically equal to zero when \(z\)~is in the domain of~\(f_{2}\). Thus \(f_{2}\)~is an odd (antisymmetric) function on the plane with a slit along the real axis from \(-1\) to~\(1\). 

Since \(f_{2}\) maps the upper half-plane bijectively onto the right-hand half-plane with a slit along the real axis from \(0\) to~\(1\) (as shown in Fig.~\ref{fig:uhp}), the reflection principle implies that \(f_{2}\) maps the whole slit plane bijectively to itself. 
If \(y\) is a positive real number, then
\begin{equation*}
-f_{2}(f_{2}(iy)) = -f_{2}( \sqrt{1+y^{2}}\,) = -(-i)\sqrt{y^{2}} =iy.
\end{equation*}
The identity principle now implies that the composite function \(-f_{2} \circ f_{2}\) is equal to the identity function. In other words, the function~\(-f_{2}\) is the inverse of~\(f_{2}\). Thus \(f_{2}\)~solves the problem of finding a holomorphic self-mapping of the slit plane with inverse function equal to its negative. 

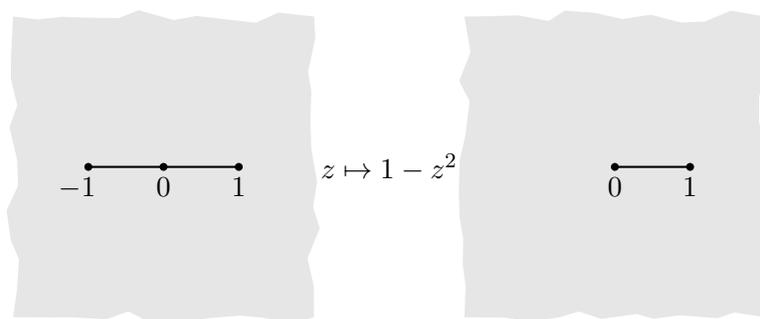
\begin{figure}
\begin{center}
\small
\begin{tikzpicture}[thick]
\fill[gray!20, decoration=random steps, decorate] (-2,-2) rectangle (2,2);
\fill (0,0) circle (1.5pt);
\fill (1,0) circle (1.5pt);
\fill (-1,0) circle (1.5pt);
\draw (3,0) node {$z\mapsto 1-z^{2}$};
\draw (-1,0) -- (1,0);
\draw (0,0) node[anchor=north] {$0$};
\draw (1,0) node[anchor=north] {$1$};
\draw (-1,0) node[anchor=north] {$-1\phantom{-}$};
\begin{scope}[xshift=6cm]
\fill[gray!20, decoration=random steps, decorate] (-2,-2) rectangle (2,2);
\draw (0,0) -- (1,0); 
\fill (1,0) circle (1.5pt);
\draw (1,0) node[anchor=north] {$1$};
\fill (0,0) circle (1.5pt);
\draw (0,0) node[anchor=north] {$0$};
\end{scope}
\end{tikzpicture}
\end{center}
\caption{An even function without a square root}
\label{fig:without}
\end{figure}

The preceding discussion demonstrates that the mocposite function \(\sqrt{1-z^{2}}\) cannot be understood as a composite function on the plane slit along the real axis from \(-1\) to~\(1\), for the function is odd instead of even. Another way to see that \(\sqrt{1-z^{2}}\) cannot be a composite function on the indicated domain is to observe that the function sending~\(z\) to \(1-z^{2}\) maps the plane slit along the real axis from \(-1\) to~\(1\) onto the plane slit along the real axis from \(0\) to~\(1\), as shown in Fig.~\ref{fig:without}. There is no holomorphic (nor even continuous) square-root function on the latter region, for the region contains the circle centered at~\(0\) with radius~\(2\), and this simple closed curve has winding number about the origin equal to~\(1\). 

On a different domain, however, the expression \(\sqrt{1-z^{2}}\) can be understood as an even composite function. Going back to~\(f_{1}\) defined on the upper half-plane, observe that \(f_{1}\)~extends continuously to the interval \( (-1,1)\) of the real axis, sending a real number~\(x\) between \(-1\) and~\(1\) to the positive square root \(\sqrt{1-x^{2}}\). By the Schwarz reflection principle, the function~\(f_{1}\) extends across this interval of the real axis to a holomorphic function~\(f_{3}\) defined on the plane with two slits, one along the real axis from \(1\) to~\(\infty\) and the other along the real axis from \(-1\) to~\(-\infty\). (See Fig.~\ref{fig:domains}(b).) Moreover, \(f_{3}({z}) = \overline{f_{3}(\overline{z})}\) for every~\(z\) in the domain. In particular, if \(y\)~is a positive real number and \(z=-iy\), then \(f_{3}(-iy) = \overline{ f_{3}(iy)} = f_{3}(iy)\) (again since \(f_{1}\), hence \(f_{3}\), takes real values on the top half of the imaginary axis). Therefore \(f_{3}\)~is an even function. 

The function sending~\(z\) to \(1-z^{2}\) maps the doubly slit plane onto the plane with a left-hand slit along the real axis from \(0\) to~\(-\infty\), which is precisely the domain of the principal branch of the square root (see Fig.~\ref{fig:squaring}), and \(f_{3}(z)\) is the composite function \(\sqrt{1-z^{2}}\). To an engineer, this function is the natural interpretation of the symbols \(\sqrt{1-z^{2}}\), not only because the function is composite but also because the reciprocal of this function is the analytic continuation to the doubly slit plane of the derivative of the inverse-sine function used in elementary differential calculus. 

In summary, the bartender's question does not admit a one-word answer. A reasonable but incomplete short answer is, ``It depends on the domain of the function.''

A deeper answer is, ``The question is wrong!'' The ultimate domain for \(\sqrt{1-z^{2}}\) is not a region in the plane but rather a two-sheeted Riemann surface, and on an abstract surface, the notions of even and odd lose meaning. The surface can be visualized as two copies of Fig.~\ref{fig:domains}(a) stitched together along the slit, the upper edge of the slit in either sheet being attached to the lower edge of the slit in the other sheet; crossing the slit corresponds to moving from one sheet to the other. Alternatively, joining two copies of Fig.~\ref{fig:domains}(b) results in an equivalent surface. 

\begin{figure}
\small
\begin{center}
\begin{tikzpicture}[thick]
\small
\fill (0,0) circle (1.5pt) (1,0) circle (1.5pt) (-1,0) circle (1.5pt);
\draw (1,0) node[anchor=west] {$1$};
\draw (-1,0) node[anchor=east] {$-1$};
\draw (0,-2) -- (0,2);
\foreach \x in {1,0.5,0.25,0.125,0.0625,0.03125} {
\begin{scope}[xscale=\x]
\draw (1,0) -- (1,2) -- (0.75,2) -- (0.75,0) -- (0.5,0);
\begin{scope}[scale=-1]
\draw (1,0) -- (1,2) -- (0.75,2) -- (0.75,0) -- (0.5,0);
\end{scope}
\end{scope}
}
\draw (0,-2.5) node {(a)};
\begin{scope}[xshift=5cm]
\foreach \x in {1,-1} {
\begin{scope}[scale=\x]
\fill (1,0) circle (1.5pt);
\draw[->,>=stealth] [domain=0:11, variable=\t, samples=110]
plot ({\t r}:{1+0.1*\t});
\end{scope}
}
\draw (1,0) node[anchor=east]{$1$};
\draw (-1,0) node[anchor=west]{$-1$};
\draw (0,-2.5) node {(b)};
\end{scope}
\end{tikzpicture}
\end{center}
\caption{Two exotic slit regions}
\label{fig:exercise}
\end{figure}
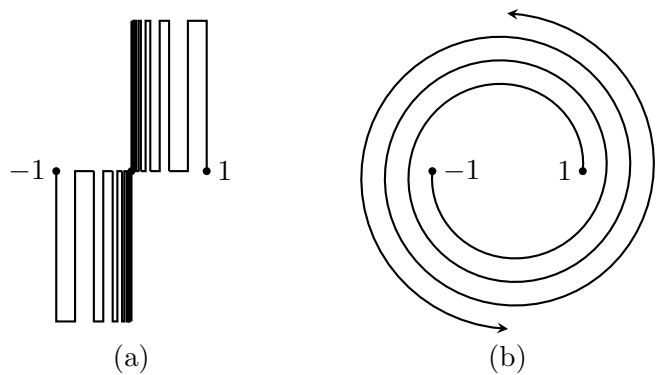

The construction cannot be implemented physically in three-dimen\-sional space, so this surface exists only in the imagination. To discuss \(\sqrt{1-z^{2}}\) with an engineer, a mathematician has to cut the Riemann surface into two pieces such that each piece projects bijectively to a planar region. The domains shown in Fig.~\ref{fig:domains} arise from two different ways of cutting the surface apart. More elaborate bisections of the surface produce exotic domains for \(\sqrt{1-z^{2}}\), such as the ones shown in Fig.~\ref{fig:exercise}. 

\begin{exercise}
On the two planar regions whose boundary slits are indicated in Fig.~\ref{fig:exercise}, is \(\sqrt{1-z^{2}}\) a mocposite function? a composite function? an even function? an odd function? 
\end{exercise}

\section{Will you join the dance?}
I invite you to seek out your own examples of mocposite functions. Such functions are easy to find; you need not travel to Wonderland \cite{alice} to encounter them. Here are a few more examples to start your feet moving the right way. 

A family of mocposite functions arises from the principal branch of the logarithm function \(\log z\), which is defined on the complex plane with a slit along the negative part of the real axis and has the property that \(e^{\log z}=z\) for every~\(z\) in the domain. The real part of \(\log z\) is equal to the natural logarithm of the modulus of~\(z\), and the imaginary part of \(\log z\) is equal to the argument (angle) of~\(z\), taken between \(-\pi \) and~\(\pi\). 

There is a holomorphic function~\(g\) on the slit plane such that \(e^{g(z)}=z^{2}\) for every~\(z\): namely, \(g(z)=2\log z\). Since \(g(z)\)~is a logarithm of~\(z^{2}\), a natural name for~\(g(z)\) is \(\log z^{2}\). This name represents a mocposite function, for \(\log z^{2}\) cannot mean the composition of a logarithm function with the squaring function. One reason is that the squaring function maps the slit plane onto the plane with a puncture at~\(0\), and there is no holomorphic logarithm function defined on the punctured plane (just as there is no holomorphic square-root function on the punctured plane). A more forceful reason is that \(g(\frac{1+i}{\sqrt{2}}) = \frac{\pi i}{2}\), but \(g(-\frac{1+i}{\sqrt{2}}) = \frac{-3\pi i}{2}\), so \(g\)~lacks the symmetry property that every function of~\(z^{2}\) must have. In particular, if \(z=-\frac{1+i}{\sqrt{2}}\), then \(\log z^{2} \ne \log(z^{2})\); ouch!

There is an analogous mocposite function \(\log z^{n}\) on the slit plane for every integer~\(n\) greater than~\(1\). More generally, a basic theorem in complex analysis says that if \(f\)~is a zero-free holomorphic function on a simply connected region of the plane (that is, a region without holes), then there exists a holomorphic function~\(g\) such that \(e^{g(z)}=f(z)\) for every point~\(z\) in the region. 

The standard proof fixes a base point~\(z_{0}\) in the region and a complex number~\(c\) such that \(e^{c} = f(z_{0})\). Set \(g(z)\) equal to \(c+\int_{z_{0}}^{z} f'(\zeta)/f(\zeta)\,d\zeta\). By Cauchy's theorem, the integral is independent of the path joining \(z_{0}\) to~\(z\) because the region is simply connected: two different paths can be deformed into each other without changing the value of the integral. The function \(fe^{-g}\) has value~\(1\) at~\(z_{0}\), and the derivative of \(fe^{-g}\) is equal to zero by the product rule, the chain rule, and the fundamental theorem of calculus. Therefore \(f=e^{g}\).

The natural name for~\(g\), a holomorphic logarithm of~\(f\), is \(\log f\). Often \(\log f\) is a mocposite function: the symbols must not be interpreted as a composition \(\log \circ f\). 

\begin{figure}
\begin{center}
\small
\begin{tikzpicture}[thick]
\draw (0,0) -- (0,1.5) (1,0) -- (1,1.5) (-1,0) -- (-1,1.5);
\fill (0,0) node[anchor=north] {$0$} circle (1.5pt);
\fill (1,0) node[anchor=north] {$\phantom{-}\pi\phantom{-}$} circle (1.5pt);
\fill (-1,0) node[anchor=north] {$-\pi\phantom{-}$} circle (1.5pt);
\draw (1.7, 0.7) node {$\ldots$};
\draw (-1.7,0.7) node {$\ldots$}; 
\end{tikzpicture}
\end{center}
\caption{Slits for a domain of \(\log\sin z\)}
\label{fig:sine}
\end{figure}
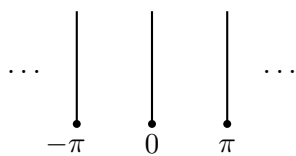

Consider, for instance, the sine function on the plane with the infinitely many unbounded vertical slits shown in Fig.~\ref{fig:sine}: for each integer~\(n\), a slit starting at the point \(n\pi\) on the real axis and going up. The zeroes of the sine function are the endpoints of the slits, so the sine function has no zero on the plane with these infinitely many slits, which is a simply connected region. Therefore a holomorphic logarithm function \(\log\sin z\) exists on the region. This function is mocposite, for the sine function maps the infinitely slit plane onto \(\mathbb{C}\setminus\{0\}\), the punctured plane, where no holomorphic logarithm function lives: composition \(\log\circ \sin\) is not defined. 

\begin{exercise}
If the value of the function \(\log\sin z\) when \(z=\frac{1}{2}\pi\) is~\(0\), then the value when \(z=2\pi+\frac{1}{2}\pi\) is~\(2\pi i\). 
More generally, if \(n\)~is an arbitrary integer, then the value of \(\log\sin z\) when \(z=(n+\frac{1}{2})\pi\) is \(n\pi i\). 
\end{exercise}

A mocposite function of a different character is the entire (holomorphic in the entire plane) function \(\cos\sqrt{z}\). This expression cannot be understood as a composite function, for \(\sqrt{z}\) is not holomorphic in a neighborhood of the origin. Nonetheless, the cosine function has a Maclaurin series containing only even powers of the variable, and replacing this variable by~\(\sqrt{z}\) produces the power series
\begin{equation*}
1-\frac{z}{2!} + \frac{z^{2}}{4!} - \frac{z^{3}}{6!} + \cdots,
\end{equation*}
which converges for every~\(z\) and thus represents an entire function that can reasonably be named \(\cos\sqrt{z}\). This function is perhaps the simplest example of an entire function of fractional order. [The order of an entire function~\(f\) is the infimum of the positive values of~\(\lambda\) for which \( f(z) e^{-|z|^\lambda}\) is a bounded function of~\(z\).] Since \(\cos z\) is the average of \(e^{iz}\) and \(e^{-iz}\), the order of \(\cos z\) evidently is~\(1\); the order of \(\cos\sqrt{z}\) is~\(1/2\). 

\section{Some hard-boiled things can be cracked}
You might think that mocposite functions are a notational curiosity of no practical importance. On the contrary, a graduate student of engineering came to me in puzzlement recently when she encountered a mocposite function in fracture mechanics. She had read in a book \cite[{\S}B.2]{fracture} about the stress intensity field induced by a crack in a material, the crack being modeled by the interval of the real axis from \(-1\) to~\(1\). The theory requires the following evaluation of an integral involving an arbitrary complex number~\(z\) lying outside the integration interval:
\begin{equation}
\frac{1}{\pi} \int_{-1}^{1} \frac{\sqrt{1-t^{2}}} {z-t}\, dt = z -  \sqrt{z^{2}-1}.
\label{crack}
\end{equation}

Since \(1-t^{2}\) is a positive real number when \(-1<t<1\), the expression \(\sqrt{1-t^{2}}\) in the integrand means the usual positive square root. Elementary real changes of variable show that the integral is an antisymmetric function of~\(z\):
\begin{equation*}
\int_{-1}^{1} \frac{\sqrt{1-t^{2}}} {-z-t}\, dt \overset{(s=-t)}{=}
\int_{-1}^{1} \frac{\sqrt{1-s^{2}}} {-z+s}\, ds
\overset{(s=t)}{=} -\int_{-1}^{1} \frac{\sqrt{1-t^{2}}} {z-t}\, dt.
\end{equation*}
Therefore the right-hand side of equation~\eqref{crack} must be antisymmetric too, but the term \(\sqrt{z^{2}-1}\) does not look antisymmetric to an engineer. As explained in \S\ref{odd}, this expression is an odd mocposite function. 

The mocposite function \(\sqrt{z^{2}-1}\) might mean either \(+i\sqrt{1-z^{2}}\) or \(-i\sqrt{1-z^{2}}\). Which choice is right for equation~\eqref{crack}? Since the integral on the left-hand side tends to~\(0\) when \(|z| \to\infty\), the expression \(\sqrt{z^{2}-1}\) needs to be close to~\(z\) when \(|z|\)~is large. The mocposite function \(\sqrt{1-z^{2}}\) constructed in~\S2 is close to \(-iz\) when \(|z|\)~is large, so \(\sqrt{z^{2}-1}\) correspondingly needs to be interpreted as \(+i\sqrt{1-z^{2}}\).

\begin{exercise}
Verify equation~\eqref{crack}, at least when \(z\) is a real number greater than~\(1\), via techniques of Calculus~II. [Suggestion: substitute \(2u/(1+u^{2})\) for~\(t\) to reduce the problem to integration of a rational function.]
\end{exercise}

\begin{figure}
\begin{center}
\small
\begin{tikzpicture}[thick]
\draw (-1,0) -- (1,0);
\draw[smooth,samples=500, decoration={markings, mark=at position 0.3 with {\arrow{<}}}, postaction={decorate}] (0,0) ellipse (1.2cm and 0.2cm);
\draw[decoration={markings, mark=at position 0.1 with {\arrow{>}}}, postaction={decorate}] (0,0) circle (2cm);
\fill (2,0) circle(1.5pt) node[anchor=west]{$R$};
\fill (0.5,1.2) circle(1.5pt) node[anchor=west]{$z$};
\end{tikzpicture}
\end{center}
\caption{An integration contour}
\label{fig:contour}
\end{figure}
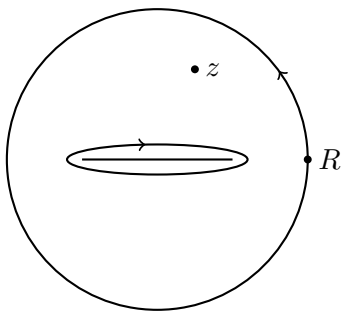

The appearance of branch issues on the right-hand side of~\eqref{crack} suggests that complex contour integration is the most natural way to evaluate the integral. One procedure is to integrate \(\sqrt{1-w^2}/(z-w)\) with respect to~\(w\) along a path consisting of a circle (oriented counterclockwise) with large radius~\(R\) and an ellipse (oriented clockwise) that surrounds the slit on the real axis (Fig.~\ref{fig:contour}). By the residue theorem, this integral equals \(-2\pi i \sqrt{1-z^2}\), or \(-2\pi \sqrt{z^2-1}\).  

On the large circle, the expression \(\sqrt{1-w^{2}}\) is \(-iw+O(1/R)\), whence the integral over the circle is 
\begin{equation*}
\int_{\text{circle}} \frac{-iw}{z-w}\,dw + O(1/R).
\end{equation*} 
By the residue theorem, the preceding expression equals \(-2\pi z+O(1/R)\). Cauchy's theorem implies that the integral over the circle is independent of~\(R\) (as long as \(R\)~is large enough that the point~\(z\) is inside the circle), so the value actually is exactly~\(-2\pi z\). 
Accordingly,
\begin{equation*}
-2\pi \sqrt{z^2-1} = -2\pi z + \int_{\text{ellipse}} \frac{\sqrt{1-w^2}}{z-w}\,d{w}.
\end{equation*}

Now let the ellipse collapse down to the slit. When \(w\)~has positive imaginary part and approaches a real value~\(t\) between \(-1\) and~\(1\), the quantity \( \sqrt{1-w^2}\) approaches the positive value \(\sqrt{1-t^2}\). The mocposite function \(\sqrt{1-w^{2}}\) is antisymmetric, so when \(w\)~has negative imaginary part and approaches a real value~\(t\) between \(-1\) and~\(1\), the quantity \( \sqrt{1-w^2}\) approaches the negative value \(-\sqrt{1-t^2}\). The top part of the ellipse approaches the slit oriented from left to right, and the bottom part of the ellipse approaches the slit oriented from right to left. Accordingly,
\begin{equation*}
 \int_{\text{ellipse}} \frac{\sqrt{1-w^2}}{z-w}\,d{w} \qquad \text{approaches}\qquad
2 \int_{-1}^1 \frac{\sqrt{1-t^2}}{z-t}\,dt.
\end{equation*}
The conclusion is that 
\begin{equation*}
-2\pi \sqrt{z^2-1} = -2\pi z+ 2 \int_{-1}^1 \frac{\sqrt{1-t^2}}{z-t}\,dt.
\end{equation*}
Dividing by \(2\pi\) shows that equation~\eqref{crack} holds. 

The trick of letting a contour collapse down to a slit when the integrand involves a (non-integer) power of \(1-t^{2}\) is an old idea. An early instance of this technique appears in the first volume of \emph{Mathematische Annalen} in a paper by Hermann Hankel containing a discussion \cite[\S3]{hankel} of integral representations of Bessel functions (special functions that appear in problems of mathematical physics involving cylindrical symmetry). One special case of Hankel's theory is the representation of the Bessel function \(J_{0}(z)\) as
\begin{equation*}
\frac{1}{\pi} \int_{-1}^{1} \frac{e^{izt}}{\sqrt{1-t^{2}}} \,dt.
\end{equation*} 
Hankel carefully explains how he understands the expression \(\sqrt{1-t^{2}}\) when \(t\)~is outside the interval of the real axis between \(-1\) and~\(1\): namely, as the product of suitably chosen branches of \(\sqrt{1-t}\) and \(\sqrt{1+t}\). A sequel to this paper~\cite{posthumous} was published two years after Hankel's untimely death at age~34 from a stroke~\cite{obituary}. Despite the clear account of branches in the original article, George Neville Watson trips up in his exposition of Hankel's work half a century later \cite[Chap.~6]{watson} by incautiously claiming non-integer powers of \(t^{2}-1\) to be even functions (and by integrating over a contour not lying in any region where \(\sqrt{t^{2}-1}\) can be defined as a holomorphic function\footnote{Experts will see how to salvage Watson's derivation by integrating a suitable holomorphic one-form over an appropriate cycle in a Riemann surface.}). 

The mocposite function on the right-hand side of equation~\eqref{crack} appears in another engineering application, one dealing with airplane wings. 
A version of the Joukowski\footnote{Famous in his native land, Nikolai Egorovich Zhukovskii (1847--1921) is ``the father of Russian aviation.'' In his French publications---notably the 1916 book \emph{A\'erodynamique}---the usual transliteration of his name is ``Joukowski,'' the spelling by which his map is commonly designated in the English literature.} 
airfoil map sends a nonzero complex number~\(z\) to the average of \(z\) and \(1/z\). At least formally, this map is an inverse of the right-hand side of equation~\eqref{crack}. Indeed, 
\begin{equation*}
 \frac{1}{z -  \sqrt{z^{2}-1}} = \frac{z + \sqrt{z^{2}-1}} {z^{2}- (z^{2}-1)} = z+\sqrt{z^{2}-1},
\end{equation*}
so, as required,
\begin{equation*}
\frac{1}{2}\left(
z -  \sqrt{z^{2}-1} + \frac{1}{z -  \sqrt{z^{2}-1}} 
\right)=z.
\end{equation*}

What is needed in addition to this formal calculation is a consideration of domains. The first observation is that the Joukowski map sending \(z\) to \(\frac{1}{2}(z+\frac{1}{z})\) is a two-to-one mapping from \(\mathbb{C}\setminus\{0\}\), the punctured plane, onto the whole plane~\(\mathbb{C}\). Indeed, if \(c\)~is an arbitrary complex number, then saying that \(\frac{1}{2}(z+\frac{1}{z})=c\) is equivalent to saying that \(z^{2}+1=2cz\), so there are two solutions for~\(z\) (counting multiplicity). Moreover, the symmetry between \(z\) and \(1/z\) reveals that the Joukowski function maps each of the regions \(\{\,z\in\mathbb{C}: 0<|z|<1\,\}\) and \(\{\,z\in\mathbb{C}: |z|>1\,\}\) one-to-one onto the same image. If \(\theta\)~is a real number, then \(\frac{1}{2}(e^{i\theta} + e^{-i\theta}) = \cos \theta\), so the Joukowski function maps the unit circle two-to-one onto the segment of the real axis between \(-1\) and~\(1\). 

Consequently, the Joukowski function maps the punctured unit disk bijectively onto the plane with a slit from \(-1\) to~\(1\) and maps the exterior of the unit disk bijectively onto the same image. The expression on the right-hand side of equation~\eqref{crack} is the inverse of one of these two functions. Since \(z -  \sqrt{z^{2}-1}\) is close to~\(0\) when the modulus of~\(z\) is large, this expression is the inverse of the restriction of the Joukowski function to the punctured unit disk. The Joukowski function is plainly odd (antisymmetric), and the inverse of an odd function is odd, so the preceding argument reconfirms the oddness of the mocposite function \(\sqrt{1-z^{2}}\). 

\section{Completed my design}
After both analysis and application, my story about mocposite functions, symmetry, and analytic continuation has come full circle. I hope that you have returned to the starting point at a new level on the Riemann surface of understanding. Here is your exit exam.

\begin{exercise}
Show that on the plane with a slit along the real axis from \(-1\) to~\(1\), the function \(\sqrt{1-\frac{1}{z^{2}}}\) is even and composite, and 
\(\sqrt{1-z^{2}}=-iz\sqrt{1-\frac{1}{z^{2}}}\). 
\end{exercise}

My secondary theme is that we mathematicians often commit expository solecisms by using confusing or ambiguous expressions, such as \(\sqrt{1-z^{2}}\), even though we purport to value rigor and precision. Lewis Carroll, from whose works I have borrowed my section titles, memorably chaffed eggheads for this shortcoming:
\begin{quote}
``When \emph{I} use a word,'' Humpty Dumpty said, in rather a scornful tone, ``it means just what I choose it to mean---neither more nor less.'' \cite{mirror}
\end{quote}
Was Humpty Dumpty a mathematician?

\end{document}